\documentclass[12pt,a4paper,reqno]{amsart}

\usepackage{amssymb}

\newcommand{\rr}{\mathbb R}
\newcommand{\zz}{\mathbb Z}
\newcommand{\nn}{\mathbb N}
\newcommand{\Tau}{\mathcal T}

\newcommand{\Per}[1]{\mathop{\mathrm{Per}_{#1}}\nolimits}
\newcommand{\Cl}[1]{\mathop{\overline{#1}}\nolimits}
\newcommand{\EssP}[1]{\mathop{{\mathcal P}(#1)}\nolimits}
\DeclareMathOperator{\Orb}{Orb}
\DeclareMathOperator{\Aper}{Aper}
\DeclareMathOperator{\lcm}{lcm}
\DeclareMathOperator{\dist}{d}

\newcounter{Stepcount}
\newenvironment{Step}{\medskip\noindent\refstepcounter{Stepcount}{\bf \theStepcount.} {}}{}
\newenvironment{LastStep}{\medskip\noindent\refstepcounter{Stepcount}{\bf \theStepcount.} {}}%
{\setcounter{Stepcount}{0}}

\theoremstyle{definition}
\newtheorem{defn}{Definition}

\theoremstyle{remark}
\newtheorem{rem}{Remark}

\theoremstyle{plain}
\newtheorem{prop}{Proposition}
\newtheorem{lemma}{Lemma}
\newtheorem{theorem}{Theorem}
\newtheorem{corr}{Corrolary}

\begin{document}

\title[One property of trajectories of Toeplitz
    flows]{One property of trajectories of Toeplitz flows}
\author[Eugene Polulyakh]{Eugene Polulyakh}
\address{Department of Topology\\
    Institute of Mathematics NAS Ukraine\\
    Tereshchenkivska str., 3\\
    Kiev\\
    Ukraine\\
    01601}
\email{polulyah@imath.kiev.ua}
\keywords{shift transform, Toeplitz flow, periodic part}
\subjclass[2000]{37B10, 37B20}

\begin{abstract}
We consider left shift transform $S$ on the space $X=\Sigma^{\zz}$ of
two-sided sequences over a compact alphabet $\Sigma$. We give an important
and sufficient condition on $x \in X$ which guarantees the~restriction of
$S$ onto orbit closure of $x$ to be a Toeplitz flow.
\end{abstract}

\maketitle

The notion of Toeplitz flow was introduced in 1969 by Jackobs and Keane in
paper~\cite{Jacobs} as certain class of subshifts of finite type. Later
this definition was expanded by S. Williams to the much more wide class of
subshifts of Bernoully shift $S$ on the space $X = \Sigma^{\zz}$ of
two-sided sequences over compact metric alphabet $\Sigma$
(see~\cite{Williams}).

Both in papers~\cite{Jacobs} and~\cite{Williams} Toeplitz flow is defined
as the restriction of $S$ onto orbit closure of so-called Toeplitz
sequence.

Let $x = (x_{n}) \in X$. Say $x_{i} \in \Sigma$ is in the {\em periodic part}
of the sequence $x$ if there exists $k \in \nn$ such that
$$
x_{i} = x_{j} \quad \text{for all} \quad j \equiv i \pmod k \,.
$$
If it is not the case we say that $x_{i}$ belongs to the {\em aperiodic part}
of the sequence $x = (x_{n})$. Sequence $x$ is called {\em Toeplitz
sequence} if it has an empty aperiodic part.

In the paper~\cite{Williams} the set of so-called {\em essential periods}
is introduced for a nonperiodic Toeplitz sequence $x$ and this set induces
in turn the {\em periodic structure} on $x$. Next, this periodic structure
defines a certain supernatural number. It appears (see~\cite{Williams})
that the flow $(\Cl{\Orb x}, S)$ admits almost one-to-one projection onto
the odometer which is defined by the same supernatural number (for the
classification of odometers by means of supernatural numbers
see~\cite{Glimm} and~\cite{Downarowicz_2}).

Toeplitz flows are remarkable since the class of all Toeplitz flows
coincides with the class of minimal flows which are symbolic and admit almost
one-to-one projection onto an odometer (for references and further
development of this result see~\cite{Downarowicz}).


The definition of Toeplitz flow is not ``homogeneous'' in the following
sense. It is known (see~\cite{Williams}) that given a Toeplitz flow $(T,
S)$ and an almost one-to-one projection $\pi : (T, S) \rightarrow G$ onto
an odometer $G$, an arbitrary point $y \in T$ is a Toeplitz sequence if
and only if $\pi$ is one-to-one in the point $y$ (i.~e. $\pi^{-1}(\pi(y)) =
\{y\}$). Hence the set of all Toeplitz sequences in $T$ is a proper
massive subset in $T$ (it contains a dense $G_{\delta}$ subset of $T$).
That's why phase space of an arbitrary Toeplits flow contains at least one
element which is not a Toeplitz sequence.

In this connection problem arises to determine whether for a given non
Toeplitz sequence $x \in X$ the dynamical system $(\Cl{\Orb x}, S)$ is a
Toeplitz flow. In the case of positive answer another problem appears
to find the periodical structure of this flow making use only of the
sequence~$x$.

We give an important and sufficient condition on $x \in X$ which guarantees
the~restriction of $S$ onto orbit closure of $x$ to be a Toeplitz flow.
Also we show how to derive the periodical structure of this flow from
$x$.

The technique applied to verify the condition allows us to expand results
of S. Williams described above (see~\cite{Williams}, theorem~2.2,
lemma~2.3 and corollary~2.4) to the case of subshifts on the space of
two-sided sequences over a Hausdorff compact alphabet (not necessarily
metrizable).

\section{Definitions and statement of results.}

Let $\Sigma$ be a  compact space, $X = \Sigma^{\zz} =
\prod_{n \in \zz} \Sigma_{n}$ with the topology of direct product.
By Tikhonov theorem $X$ is also the compact space.
We write elements of $X$ as $x = (x_{n})$.

In the case when $(\Sigma, \rho)$ is a metric space the distance
$$
d(x, y) = \sum_{n \in \zz} 2^{-|n|} \rho(x(n), y(n))
$$
is known to induce the product topology on $X$.

In what follows the next property of $X$ will be useful to us
(see.~\cite{Kelley}).
\begin{prop}\label{prop_1}
A sequence of points $\{x_{i}\}$ in a product $\prod_{n \in \zz}
X_{n}$ of topological spaces converges to $x \in \prod_{n \in \zz}
X_{n}$ if and only if the sequence $\{x_{i}(n)\}$
converges to $x(n)$ for every $n \in \zz$.
\end{prop}

Let us designate by $S: X \rightarrow X$ the left shift homeomorphism $S(x(n)) =
x(n+1)$, $n \in \zz$.

For $x \in X$, $p \in \nn$ and $\sigma \in \Sigma$ let
$$
\Per{p}(x, \sigma) = \left\{ \; n \in \zz \;|\; x(n') = \sigma \text{{} for all
} n' \equiv n \pmod p \; \right\} \,,
$$
$$
\Per{p}(x) = \bigcup_{\sigma \in \Sigma} \Per{p}(x, \sigma) \,,
$$
$$
\Aper(X) = \zz \setminus \Bigl( \bigcup_{p \in \nn} \Per{p}(x) \Bigr) \,.
$$
By $p$-skeleton of $x$ we shall name that part of a sequence $(x(n))$, which
has the period $p$.

Let us designate
$$
M_{p}(x) = \max \{ \, k \in \nn \;|\; \exists n \in \zz: n+i \in
\Per{p}(x), \; i = 0, 1, \ldots, k-1 \, \} \,.
$$
In other words $M_{p}(x)$ is the maximal length of a block contained in
$p$-skeleton of $x$. Note, that $M_{p}(x) = \infty$ for periodic
sequence $x$ with period $p$ and $M_{p}(x) < p$, if the
sequence $x$ is not periodic.

Let us remind some important definitions.

\begin {defn}
Sequence $\eta \in X$ is called {\em Toeplitz}, if $\Aper(\eta) = \emptyset$
(in this case the dynamic system $(\Cl{\Orb(\eta)}, S)$ is also referred
as {\em Toeplitz}).
\end {defn}

\begin{defn}
Let $(X, F)$ be a dynamic system with discrete time, $x \in X$.
The point $x$ is {\em recurrent}, if for its arbitrary open
neighbourhood $U$ there exists $n(U)$, such that for any $k \in \zz$
$$
U \cap \Bigl( \bigcup_{i = k}^{k+n(U)-1} \{ F^{i} (x) \} \Bigr) \neq
\emptyset \,.
$$
\end{defn}

\begin{defn}
Let $(X, F)$ be a dynamic system with discrete time, $x \in X$.
The point $x$ is said to be {\em almost periodic}, if for its arbitrary open
neighbourhood $U$ we can find $n(U) \in \nn$, such that
$$
\bigcup_{k \in \zz} \{F^{k n(U)} (x) \} \subset U \,.
$$
\end{defn}

Clearly each periodic sequence $x \in X $ is
Toeplitz. It is easy to check that every Toeplitz sequence is almost periodic
since each block of
Toeplitz sequence $\eta$ is contained in its $p$-skeleton for some
$p$. Hence, according to Birkgoff theorem $\Cl{\Orb(\eta)}$
is a minimal set of dynamic system $(X, S)$ (see~\cite{Alekseev,Morse}).

Let $\eta \in X$ be an aperiodic Toeplitz sequence.
Generally speaking an equality $\Aper(x) = \emptyset$
is not carried out for an arbitrary $x \in \Cl{\Orb(\eta)}$.

Consider a special case when $\Sigma$ is metric space. From one hand,
every Toeplitz flow $(\Cl{\Orb(\eta)}, S)$ in $X$ have to be expansive
(see~\cite{G-J}). From the other hand, every odometer is an equicontinuous
dynamic system (see remark~\ref{rem_odom} below). It is known that any
Toeplitz flow admits almost one-to-one projection onto odometer and such
a projection have to be one-to-one precisely in points which are Toeplitz
sequences (see~\cite{Williams}). So, if every point of a certain Toeplitz
flow is a Toeplitz sequence, then this Toeplitz flow have to be conjugate
to an odometer. In particular it must be expansive and equicontinuous
simultaneously, and this is impossible.

Properties of sequences from $\Cl{\Orb(\eta)}$ are in details investigated
in~\cite{Williams}. However it is not known, what should be properties of
the point $x \in X$ the set $\Cl{\Orb(x)}$ to contain some Toeplitz
sequence.

The answer to this question gives the following
\begin{prop} \label{prop_2}
If for some $x \in X$
\begin{equation} \label{eq_1}
\limsup_{p \rightarrow \infty} M_{p} (x) = \infty \;,
\end{equation}
then set $\Cl{\Orb(x)}$ contains a Toeplitz sequence.

If in addition point $x$ is recurrent, then dynamic system
$(\Cl{\Orb(x)}, S)$ is Toeplitz.
\end {prop}

There is a natural question: what additional information it is possible
to take about structure of dynamic system $(\Cl{\Orb(x)}, S)$ under condition
that a point $x$ is recurrent?

Let's remind definition of periodic structure of a Toeplitz sequence
(see~\cite{Williams}).

\begin{rem}
Let $x \in X$. If $p \,|\, q$, then $\Per{p}(x) \subseteq \Per{q}(x)$.
\end{rem}

\begin{defn}
Let's call $p$ {\em the essential period} of a sequence $x$, $p \in
\EssP{x}$, if for any $q \in \nn$
$$
(\Per{p}(x, \sigma) \subset \Per{p}(x, \sigma) -
q \quad \forall \sigma \in \Sigma)
\quad \Rightarrow \quad
( p \, | \, q) \;.
$$
\end{defn}
In other words, $p \in \EssP{x}$ if and only if $p$-skeleton of $x$
is not periodic for any smaller period.

\begin{rem} \label{rem_2}
It is easily checked, that if $p$, $q \in \EssP{x}$, then $\lcm (p, q)
\in \EssP{x}$ (see~\cite{Williams}).
\end{rem}

\begin{defn} \label{per_struct_1}
{\em Periodic structure} of a nonperiodic Toeplitz
sequence $\eta$ is the growing sequence
$\{p_{i}\}_{i \in \nn}$ of natural numbers, such that
\begin{itemize}
    \item[(i)] $p_{i} \in \EssP{\eta}$ for all $i \in \nn $;
    \item[(ii)] $p_{i} \,|\, p_{i+1}$;
    \item[(iii)] $\bigcup_{i \in \nn} \Per{p_{i}}(x) = \zz$.
\end{itemize}
\end{defn}
To within the equivalence relation which we will not describe here
periodic structure for a Toeplitz sequence
is determined uniquely (see~\cite{Williams}). For our purposes it is enough to know
that any subsequence of a sequence from previous
definition sets an equivalent periodic structure.

Now we shall determine periodic structure for any recurrent
sequence $x \in X$ which satisfies
the relation~(\ref{eq_1}).

\begin{rem} \label{rem_3}
Let $p \,|\, q$ for some $p$, $q \in \nn$. Then $M_{p}(x) \leq
M_{q}(x)$.
\end{rem}

\begin{defn} \label{per_struct_2}
{\em Periodic structure} of an aperiodic
sequence $x$ which satisfies
the relation~(\ref{eq_1}) is the growing sequence
$\{p_{i}\}_{i \in \nn}$ of natural numbers, such that
\begin{itemize}
    \item[(i)] $p_{i} \in \EssP{\eta}$ for all $i \in \nn$;
    \item[(ii)] $p_{i} \,|\, p_{i+1}$;
    \item[(iii$'$)] $\lim_{i \rightarrow \infty} M_{p_{i}}(x) = \infty$.
\end{itemize}
\end{defn}

\begin{prop} \label{prop_3}
For each aperiodic sequence $x$, which
satisfies to the relation~(\ref{eq_1}), there exists some periodic
structure.
\end {prop}

For the benefit of such definition of periodic structure speak the following
results.

\begin{theorem} \label{th_1}
Suppose the sequence $\{p_{i}\}$ determines certain periodic
structure (in the sense of definition~\ref{per_struct_2}) for a recurrent
sequence $x \in X$ which satisfies relation~(\ref{eq_1}). Then
there exists a Toeplitz sequence $\eta \in X$, such that
$\Cl{\Orb(x)} = \Cl{\Orb(\eta)}$ and the sequence $\{p_{i}\}$ evaluates
periodic structure on $\eta$ (in sense of definition~\ref{per_struct_1}).
\end{theorem}

\begin{corr}
Let $x \in X $ be a recurrent sequence satisfying
to equality~(\ref{eq_1}). Then the periodic structure for $x$ is determined
uniquely (to within the relation of equivalence from~\cite{Williams}).
\end{corr}

\section {Proof of the main results.}

\subsection {Proof of proposition~\ref{prop_3}.}

We fix $x \in X$. Divide the proof into several steps.

\begin{Step}
Suppose $\Per{p}(x) \neq \emptyset$ for some $p \in \nn$. Find
minimal $k \in \nn $, such that $k \,|\, p$ and $ \Per{p}(x) =
\Per{k}(x)$.

Let us check that $k \in \EssP{x}$. Two lemmas will be necessary
for this purpose.

\begin{lemma}
Suppose that the following condition
$$
\Per{p}(x, \sigma) = \Per{p}(x, \sigma) + m_{i} \quad \forall \sigma \in
\Sigma \,, \; i = 1, 2
$$
is satisfied for some $m_{1}$, $m_{2} \in \nn$.
Let $b \in \nn$, $0 \leq b \leq m_{2}-1$, be a remainder of the division of $m_{1}$
into $m_{2}$. Then
$$
\Per{p} (x, \sigma) = \Per{p} (x, \sigma) + b
\quad \forall \sigma \in \Sigma \,.
$$
\end{lemma}

\begin{proof} On a condition $m_{1} = a m_{2} + b$, $a \in \zz_{+}$.
For every $\sigma \in \Sigma$
$$
\Per{p} (x, \sigma) = \Per{p} (x, \sigma) + m_{1} =
\left(\Per{p} (x, \sigma) + a m_{2} \right) + b =
\Per{p} (x, \sigma) + b.
$$
\end{proof}

\begin{lemma} \label{lemma_2}
Let for some $q \in \nn$ following condition is satisfied
$$
\Per{p} (x, \sigma) = \Per{p} (x, \sigma) +
q\quad \forall \sigma \in \Sigma \,.
$$
Then $\Per{\gcd (p, q)}(x) = \Per{p}(x)$.
\end{lemma}

\begin{proof} Consider Euclidean algorithm of a finding of $\gcd (p, q)$.
\begin{equation} \label{eq_2}
\begin{split}
    p & = a_{1} q + b_{1} \,, \;\:\qquad 0 \leq b_{1} < q \,; \\
    q & = a_{2} b_{1} + b_{2} \,, \,\qquad 0 \leq b_{2} < b_{1} \,; \\
      & \cdots \\
    a_{n-2} & = a_{n} b_{n-1} + b_{n} \,, \quad 0 \leq b_{n} =
        \gcd (p, q) < b_{n-1} \,; \\
    a_{n-1} & = a_{n+1} b_{n} \,.
\end{split}
\end{equation}
Applying the previous lemma by turns to each line of~(\ref {eq_2})
we are convinced that for $i = 1, \ldots, n$
$$
\Per{p} (x, \sigma) = \Per{p} (x, \sigma) + b_{i} \quad \forall \sigma \in \Sigma \,.
$$
In particular, $\Per{p} (x, \sigma) = \Per{p} (x, \sigma) + \gcd (p, q)$
$\forall \sigma \in \Sigma$.

Hence, $\Per{p} (x) = \bigcup_{\sigma \in \Sigma} \Per{p} (x, \sigma)
\subseteq \Per{\gcd (p, q)}(x)$. On the other hand,
since $\gcd (p, q) \,|\, p$ the opposite inclusion
$\Per{\gcd (p, q)}(x) \subseteq
\Per{p}(x)$ is also true.
\end{proof}

So, let $k$ be the minimal from divisors of $p$, such that $\Per{k}(x) =
\Per{p}(x)$. Let for some $q \in \nn$ the equality
$$
\Per{k} (x, \sigma) = \Per{k} (x, \sigma) + q \quad \forall \sigma \in \Sigma
$$
is hold true.
Then $\Per{\gcd (k, q)} (x) = \Per{k}(x) = \Per{p}(x) $ on
lemma~\ref {lemma_2}. Since $k \,|\, p$ and $\gcd (k, q) \,|\,
k$ then $\gcd (k, q) = k$ by virtue of a choice of $k$ and
$k \,|\, q$. That is $k \in \EssP{x}$.

\begin{rem} \label{rem_4}
As $ \Per{k}(x) = \Per{p}(x)$  on the construction then $M_{k}(x) = M_{p}(x)$.
\end{rem}
\end{Step}

\begin{LastStep}
Now we shall proceed directly to the construction of the periodic structure for
$x$.

Taking into account the equality~(\ref{eq_1}) we shall choose a sequence
$\{p_{i}\}_{i \in \nn}$ of the natural numbers to comply with the relation
$$
\lim_{i \rightarrow \infty} M_{p_{i}}(x) = \infty \,.
$$

Further, using argument stated before we shall choose
the least divisor $k_{i}$ of $p_{i}$ for every $i \in \nn$ such that $\Per{k_{i}} (x) =
\Per{p_{i}}(x)$. We shall receive a sequence $\{k_{i}\}$ of the essential
periods for $x$ satisfying the relation $\lim_{i \rightarrow \infty} M_{k_{i}}(x) =
\infty$ (see remark~\ref {rem_4}).

Set
$$
q_{i} = \lcm (k_{1}, k_{2}, \ldots, k_{i})
$$
for every $i \in \nn $.
It is easily verified that $q_{i} \,|\, q_{i+1}$, $i \in \nn$.
Remark~\ref{rem_2} guarantees that a sequence $\{q_{i}\}$
contains only the essential periods for $x$, and the equality
$$
\lim_{i \rightarrow \infty} M_{q_{i}}(x) = \infty
$$
follows from remark~\ref{rem_3}.

Proposition~\ref {prop_3} is completely proved.
\end{LastStep}

\subsection {Proof of theorem~\ref{th_1} and proposition~\ref{prop_2}.}

We fix periodic structure $\{q_{i}\}$ on $x$. Passing to
subsequence it is possible to suppose that
\begin{equation} \label{eq_3}
M_{q_{i+1}}(x) \geq 3 q_{i} + M_{q_{i}}(x) \,, \quad i \in \nn \,.
\end{equation}

First we shall construct a Toeplitz sequence $\eta \in \Cl{\Orb (x)}$
such that
$$
\zz = \bigcup_{i \in \nn} \Per{q_{i}}(\eta)
$$
( and thus we shall prove proposition~\ref {prop_2}), and then we shall show that
$q_{i} \in \EssP{\eta}$, $i \in \nn$.

\begin{Step}
We fix a sequence $\{m_{i}\}_{i \in \nn}$ of integers, such that $m_{i} +
j \in \Per{q_{i}}(x)$ for all $i \in \nn$ and $j \in \{0, 1, \ldots,
M_{q_{i}}(x) - 1\} $, that is for every $i \in \nn$ if $n \equiv 0 \pmod
{q_{i}}$ then
$$
x(m_{i} +j) = x(m_{i}+j+n) \,, \quad
j = 0, 1, \ldots, M_{q_{i}}(x) - 1 \,.
$$

From the relation~(\ref{eq_3}) it follows that
$$
\left[ M_{q_{i+1}}(x) - (q_{i} + M_{q_{i}}(x)) \right] - q_{i} \geq q_{i}
\,,
$$
therefore for every $i \in \nn$ there exists
$$
s_{i} \in \left[ m_{i+1} + q_{i} \,,\; (m_{i+1} + M_{q_{i+1}}(x)) - (q_{i} + M_{q_{i}}(x))
\right] \,,
$$
which complies with the equality $m_{i} \equiv s_{i} \pmod {q_{i}}$.

Let us designate
$$
d_{l}(i) = s_{i} - m_{i+1} \,, \quad
d_{r}(i) = \left( m_{i+1} + M_{q_{i+1}}(x) \right) - (s_{i} +
M_{q_{i}}(x)) \,.
$$
Note that $d_{l}(i)$ and $d_{r}(i)$ are the numbers of symbols of the block
$$
x(m_{i+1}), x(m_{i+1}+1), \ldots, x(m_{i+1} + M_{q_{i+1}}(x) - 1) \,,
$$
standing accordingly at the left and at the right of the block
$$
x(s_{i}), x(s_{i}+1), \ldots, x(s_{i} + M_{q_{i}}(x) - 1)
$$
of the sequence $x = (x(n))$.

It is not difficult to see that
\begin{equation}
\begin{split}
    d_{l}(i) & \geq q_{i} \,, \\
    d_{r}(i) & \geq
        (m_{i+1} + M_{q_{i+1}}(x)) - \\
        & \quad - \left[(m_{i+1} + M_{q_{i+1}}(x)) -
        (q_{i} + M_{q_{i}}(x)) + M_{q_{i}}(x) \right] = q_{i} \,.\\
\end{split}
\end{equation}

Consider a sequence of integers
\begin{align*}
    k_{1} & =  m_{1} \,, \\
    k_{2} & =  k_{1} + (s_{1} - m_{1}) = s_{1} \,, \\
    &  \cdots \\
    k_{j} & =  k_{j-1} + (s_{j-1} - m_{j-1}) = \\
        & =  m_{1} + (s_{1} - m_{1}) + \ldots + (s_{j-1} - m_ {j-1}) \,,
        \quad j > 1 \,.\\
\end{align*}
and a sequence $z_{j} = S^{k_{j}}(x)$, $j \in \nn$, of elements of the set
$\Cl{\Orb{x}}$. We shall note obvious equalities
\begin{gather*}
    z_{j} = S^{k_{j} - k_{l}} \circ S^{k_{l}}(x) = S^{k_{j} -
        k_{l}}(z_{l}) \,, \\
    k_{j} - k_{l} = \sum_{i = l}^{j-1}(s_{i} - m_{i}) \,, \quad l < j
        \,.
\end{gather*}

Notice that since for all $j \in \nn$ on construction
$q_{j} \,|\, (s_{j}-m_{j})$ and $q_{j} \,|\, q_{j+1}$ then
$$
q_{l} \,|\, (k_{j}-k_{l}) \,, \quad l < j
$$
and for every $n \in \Per{q_{l}}(x)$ and $j > l$ we have
$$
z_{j}(n) = z_{l}(n) = x(n+k_{l}) \,.
$$

Let us designate
\begin{equation}
P_{l} = \Per{q_{l}}(z_{l}) = \Per{q_{l}}(x) - k_{l} \,.
\end{equation}

Above we have already checked up that $z_{j}(n) = z_{l}(n)$ for all $j > l$ and $n \in
P_{l}$. We shall show now that $\bigcup_{l \in \nn} P_{l} = \zz$.

On construction $[m_{l}, m_{l} + M_{q_{l}}(x)-1] \subset \Per{q_{l}}(x)$, $l \in
\nn$, hence $[m_{l} - k_{l}, m_{l} + M_{q_{l}}(x) - 1 - k_{l}] \subset
P_{l}$.

Notice that for $l = 1$
\begin{align*}
    m_{1} - k_{1} &= 0 \,, \\
    m_{1} + M_{q_{1}}(x)-1-k_{1} &= M_{q_{1}}(x) - 1 \,,
\end{align*}
hence $[0, M_{q_{1}}(x) - 1] \subset P_{1}$.

When $l \geq 2$ we have
\begin{equation*}
\begin{split}
    m_{l}-k_{l} & = -m_{1} - (s_{1}-m_{1}) - \ldots - (s_{l-1}-m_{l-1})
        + m_{l} = \\
    & = - (s_{1}-m_{2}) - \ldots - (s_{l-1}-m_{l}) = \\
    & = -d_{l}(1) - \ldots - d_{l}(l-1) \leq \\
    & \leq -q_{1} - \ldots - q_{l-1} \,;
\end{split}
\end{equation*}
\begin{multline*}
m_{l} + M_{q_{l}}(x) - 1 - k_{l} = \\
\begin{split}
& = m_{l} + M_{q_{l}}(x) - 1 - m_{1} -
    \sum_{i=1}^{l-1} (s_{i} - m_{i}) = \\
& = M_{q_{l}}(x) - 1 + \sum_{i=2}^{l} (m_{i} - s_{i-1}) = \\
& = M_{q_{1}}(x) - 1 + \sum_{i=2}^{l}
    \left[ (m_{i}+M_{q_{i}}(x)) - (s_{i-1}+M_{q_{i-1}}(x)) \right] = \\
& = M_{q_{1}}(x) - 1 + \sum_{i=1}^{l-1} d_{r}(i) \geq
    M_{q_{1}}(x) - 1 + \sum_{i=1}^{l-1} q_{i}\,.
\end{split}
\end{multline*}
Hence, for all $l > 1$
\begin{equation}
\left[ -\sum_{i=1}^{l-1} q_{i} \,, \; M_{q_{1}}(x) - 1 + \sum_{i=1}^{l-1}
q_{i} \right] \subset P_{l} \,.
\end{equation}
On construction $q_{i} \geq 1$, $i \in \nn$, so $\bigcup_{l \in \nn} P_{l} =
\zz$.

Therefore, it is correctly determined $\eta \in X$ which meets the equality
$\eta(n) = z_{l}(n)$, if $n \in P_{l}$.

It is easy to see that $P_{l} \subseteq \Per{q_{l}}(\eta)$, $l \in \nn $.
Furthermore, from proposition~\ref {prop_1} it follows that
$\eta = \lim_{i \rightarrow \infty} z_{i}$.

\begin{rem}
So, we have constructed the Toeplitz sequence $\eta \in \Cl{\Orb {x}}$.
In the argument above we have nowhere used recurrence of $x$.
\end{rem}

Let now a point $x$ is recurrent. Under Birkgoff theorem the set
$\Cl{\Orb{x}}$ is minimal, hence $\Cl{\Orb{\eta}} = \Cl{\Orb{x}}$.
\end{Step}

\begin{LastStep}
Let $\eta \in X$ be a Toeplitz sequence, $x, y \in
\Cl{\Orb{\eta}}$. Let us prove that $M_{p}(x) = M_{p}(y)$ for every $p \in
\nn$ and $\EssP{x} = \EssP{y}$.

\begin{lemma} \label{lemma_3}
Let $A$ be a minimal subset of dynamic system $ (X,
S)$, $x, y \in A$. Let $\Per{p}(x) \neq \emptyset$ for some $p \in
\nn$. Then there exists $n(p) \in \zz$ which satisfies the conditions
\begin{itemize}
    \item[(i)] $\Per{p}(y) = \Per{p}(x) - n(p)$;
    \item[(ii)] $x(k+n(p)) = y(k)$ for every $n \in \Per{p}(x)$;
\end{itemize}
\end{lemma}

\begin{proof}
{\bf 1.} First we shall prove that there exists $n(p) \in \zz$
which satisfies to a condition (ii) (hence for this $n(p)$
inclusion $\Per{p}(x) \subseteq \Per{p}(y) + n(p)$ is valid).

Since the set $A$ is minimal then $A = \Cl{\Orb{x}} = \Cl{\Orb{y}}$ and
there exists a sequence $\{z_{j} = S^{k_{j}} (x)\}_{j \in \nn}$
converging to a point $y$.

Let us say that $k_{i} \sim k_{j}$ if $k_{i} \equiv k_{j} \pmod p$.
Under this relation the set $\{k_{i}\}$ will fall into no more than on $p$ classes
of equivalence. Obviously, at least one of these classes contains
infinite number of elements. Hence, passing to a
subsequence we can assume that $k_{i} \equiv k_{j} \pmod p$
for all $i, j \in \nn$.

Then $ \Per{p}(z_{i}) = \Per{p}(z_{1}) = \Per{p}(x) - k_{1}$ for all $i
\in \nn$ (we shall designate $P(p) = \Per{p}(x) - k_{1}$). Moreover, $z_{i}(k) =
z_{1}(k) = x(k+k_{1})$ for all $k \in P(p)$.

From proposition~\ref {prop_1} it follows that $y(k) = z_{1}(k) = x(k+k_{1})$
for every $k \in P(p)$. Therefore, $(\Per{p}(x) - k_{1}) \subseteq
\Per{p}(y)$ also it is possible to let $n(p) = k_{1}$.

{\bf 2.} Let us check now a correlation~(i).

Assume that $\Per{p}(x) \subsetneqq \Per{p}(y) + n(p)$. Repeating the
argument of item~1 and changing roles of $x$ and $y$, we shall find $m(p) \in
\zz$ such that $\Per{p}(y) \subseteq \Per{p}(x) + m(p)$. Then
$\Per{p}(y) + n(p) \subseteq \Per{p}(x) + (m(p) + n(p))$ and
$$
\Per{p}(x) \subsetneqq \Per{p}(x) + (m (p) + n (p)) \,.
$$
Clearly, $m(p) + n(p) \neq 0$.

Obviously, for every $r \in \zz$
\begin{equation} \label{eq_7}
\Per{p}(x) + r \; \subsetneqq \; \Per{p}(x) + (m(p) + n(p)) + r \,.
\end{equation}

Let $s = \lcm (m(p) + n(p), p)$. Then $s = a(m(p) + n (p))$ for
certain $a \in \zz \setminus \{ 0 \}$.

Assume that $a < 0$ (the case $a > 0$ is examined
similarly). Using a relation~\ref{eq_7}, we shall receive the following chain
of inclusions
\begin{align*}
\Per{p}(x) & \supsetneqq \Per{p}(x) - (m(p) + n(p)) \supsetneqq \\
& \supsetneqq \Per{p}(x) - 2(m(p) +n(p)) \supsetneqq \\
& \ldots \\
& \supsetneqq \Per{p}(x) + a(m(p) + n(p)) \,.
\end{align*}
However, on construction $p \,|\, s$, hence
$$
\Per{p}(x) + a(m(p) + n(p)) = \Per{p}(x) + s = \Per{p}(x)
$$
by definition of $\Per{p}(x)$.

The received contradiction finishes the proof of lemma.
\end{proof}

\begin{corr} \label{corr_2}
$M_{p}(x) = M_{p}(y)$ for every $p \in \nn$ and $\EssP{x} = \EssP{y}$.
\end{corr}

Applying now lemma~\ref{lemma_3} and corollary~\ref{corr_2} to the
sequence $\{q_{i}\} $ we verify that $P_{i} =
\Per{q_{i}}(\eta)$ and $q_{i} \in \EssP{\eta}$, $i \in \nn$. For the completion of
the proof of theorem~\ref{th_1} it remains to recall
the equality $\bigcup_{i \in \nn} P_{i} = \zz$ which we have already checked
above.
\end{LastStep}
\qed

\section{Toeplitz subshifts on the space of two-sided sequences over a
Hausdorff compact alphabet}

\subsection{Odometers and periodic partitions of dynamic systems}

\begin{defn}
A non-bounded sequence $\{a_{i} \in \nn\}_{i \in \nn}$ is called
{\em regular} if $a_{i}$ divides $a_{i+1}$ for every $i \in \nn$.
\end{defn}

We fix regular sequence
$\{n_{i} \in \nn\}_{i \in \nn}$ (without loss of generality
it is possible to assume that $n_{i+1} \neq n_{i}$, $i \in \nn$).

Let us consider a sequence of finite cyclic groups $\zz_{n_{i}} =
\zz / n_{i} \zz$ and group homomorphisms
$$
\varphi_{i} : \zz_{n_{i+1}} \rightarrow \zz_{n_{i}} \;,
$$
$$
\varphi_{i} : 1 \mapsto 1 \;.
$$
Let us take an inverse limit $A = \projlim_{i \rightarrow \infty}
\zz_{n_{i}}$ of this sequence of groups and homomorphisms.
We receive an abelian group $(A, +)$.

Provide each set $\zz_{n_{i}} = \{0, 1, \ldots, n_{i} -1 \}$
with discrete topology. Each of maps $\varphi_{i}$ is
continuous in this topology. Space $A$ with topology $\Tau$
of the inverse limit is homeomorphic to a Cantor set $\Gamma$.

It is easy to see, that in the group $(A, +)$ operation of addition and pass to
an opposite element are continuous in the topology $\Tau$, thus $A$
turns to be a continuous group.

\begin{rem}
We remind that an inverse limit $A = \projlim_{i \rightarrow \infty}
\zz_{n_{i}}$ could be imagined as a subset
\begin{equation} \label{eq_3.1}
A = \{\vec {a} = (a_{i} \in \zz_{n_{i}}) \; | \; \varphi_{i} (a_{i+1}) = a_{i}
\,, \; i \in \nn \}
\end{equation}
of the direct product
\begin{equation} \label{eq_3.2}
\prod_{i \in \nn} \zz_{n_{i}} \;.
\end{equation}

In such notation the operation of addition in $A$ is defined component-wise, that is
$\vec{a} + \vec{b} = (a_{i} + b_{i})$ for any $\vec{a} = (a_{i})$,
$\vec{b} = (b_{i}) \in A$.
\end{rem}

As is known, the topology of the direct product~(\ref{eq_3.2}) is set
through a basis consisting of so-called cylindrical sets
\begin{gather*}
U(x_{i_{1}}, \ldots, x_{i_{k}}) =
\{(a_{i}) \; | \; a_{i_{s}} = x_{i_{s}} \,, \; s = 1, \ldots, k \} \;; \\
x_{i_{s}} \in \zz_{n_{i_{s}}} \,, \; i_{1} < \ldots < i_{k} \,, \; k \in \nn \;.
\end{gather*}

From definition of the set $A$ (see relation~(\ref{eq_3.1})) it is easy
to see that
$$
U(x_{i_{1}}, \ldots, x_{i_{k}}) \cap A = U (x_{i_{k}}) \cap A
$$
for any $k \in \nn $, $i_{1} < \ldots < i_{k}$ and $x_{i_{s}} \in
\zz_{n_{i_{s}}}$.
So, the family of sets
\begin{equation} \label{eq_3.3}
\begin{split}
    V_{x_{j}} & = U(x_{j}) \cap A = \{(a_{i}) \in A \; | \;
        a_{j} = x_{j} \} \quad = \\
    & = \{(a_{i}) \in A \; | \; a_{j} = x_{j} \,, \; a_{k} = \varphi_{k} \circ
        \ldots \circ \varphi_{j-1} (x_{j}) \mbox{ when } k < j \} \;; \\
    & \qquad j \in \nn \,, \; x_{j} \in \zz_{n_{j}}
\end{split}
\end{equation}
is base of the topology of space $A$.

The {\em natural metric} $\dist : A \times A \rightarrow \rr_{+}$ on $A$
associated with the sequence $\{n_{i}\}$ is defined as follows
$$
\dist(\vec{x}, \vec{y}) = \frac{1}{m} \;, \quad
m = \min \{i \in \nn \; | \; x_{k} = y_{k} \mbox{ when } k < i \mbox{ and }
x_{i} \neq y_{i} \} \;.
$$
The correctness of this definition is checked immediately.

Consider an element $\vec{e} = (1) = (1, \ldots, 1, \ldots) \in A$. This
element is called {\em generator} of group $A$ and has the property
that the cyclic subgroup $\langle \vec{e} \rangle$ generated by it is dense
in $A$ in the topology $\Tau$.

Obviously, shift mapping
$$
g: A \rightarrow A \;, \qquad
g: \vec{x} \mapsto \vec{x} + \vec{e} \;,
$$
is a homeomorphism.

\begin{defn}
Dynamic system $(A, g)$ is called an {\em odometer}.
\end{defn}

\begin{rem}
From the fact the subgroup $\langle \vec{e} \rangle$ is dense in $A$ it immediately
follows that each trajectory of d. s. $(A, g)$ is dense in $A$, that is
odometer always is a minimal dynamic system.
\end{rem}

\begin{rem} \label{rem_odom}
It is easy to verify that in the natural metric defined above the mapping
$g$ is isometric. Specially, the family of mappings $\{g^{k}\}_{k \in
\zz}$ is equicontinuous, so the odometer $(A, g)$ is the equicontinuous
dynamic system.

Actually, it is known that odometers are precisely all equicontinuous
minimal dynamic systems on the Cantor set.
\end{rem}

Assume a compact Hausdorff space $X$ and homeomorphism $f: X \rightarrow
X$ are given.

\begin{defn} \label{defn_9}
We call a finite family $\{ W_{i} \}_{i=0}^{n-1}$ of subsets of space
$X$ a {\em periodic partition} of the dynamic system $(X, f)$ {\em of
length} $m$, if it satisfies to the following requirements:
\begin{itemize}
    \item[(i)] all $W_{i}$ are open-closed subsets of $X$;
    \item[(ii)] $W_{i} = f(W_{i-1})$, $i = 1, \ldots, n-1$ and $W_{0} =
        f(W_{n-1})$;
    \item[(iii)] $W_{i} \cap W_{j} = \emptyset$ when $i \neq j$;
    \item[(iv)] $X = \bigcup_{i=0}^{n-1} W_{i}$.
\end{itemize}
\end{defn}

\begin{lemma}
Assume $(A, g)$ is an odometer built with the help of a regular sequence
$\{n_{i}\}_{i \in \nn}$.

For any $k \in \nn$ and $x_{k} \in \zz_{n_{k}}$ a family of sets
$\{W^{(n_{k})}_{j} = V_{x_{k} +j}\}_{j = 0, \ldots, n_{i}-1}$ forms
periodic partition of a dynamic system $(A, g)$ of length $n_{k}$.
\end{lemma}

\begin{proof}
Obviously,
$$
A = \bigcup_{s \in \zz_{n_{k}}} V_{s} = \bigcup_{j \in \zz_{n_{k}}}
V_{x_{k} +j} \;.
$$
Hence, for the family $\{W^{(n_{k})}_{j}\}$ the requirement~(iv)
of Definition~\ref{defn_9} is carried out.

Since all sets $V_{x_{k}+j}$, $j \in \zz_{n_{k}}$, are open on
definition and pairwise disjoint, family $\{W^{(n_{k})}_{j}\}$
satisfies also to properties~(i) and~(iii) of a periodic partition.

For completion of the proof we need to verify that $g(V_{a_{k}}) =
V_{a_{k}+1}$ (here $1 \in \zz_{n_{k}}$) for every $a_{k} \in \zz_{n_{k}}$.

Let $\vec{b} = (b_{i}) \in V_{a_{k}}$. Then $b_{k} = a_{k}$ and
$g(\vec{b}) = \vec{b} + \vec{e} = (b_{i}+1) \in V_{a_{k}+1}$.
Hence, $g(V_{a_{k}}) \subseteq V_{a_{k}+1}$.

Back, let $\vec{c} = (c_{i}) \in V_{a_{k}+1}$. Then $c_{k} =
a_{k}+1$ and $g^{-1} (\vec{b}) = \vec{c} - \vec{e} = (c_{i}-1) \in
V_{a_{k}}$. Hence, $g (V_{a_{k}}) \supseteq V_{a_{k}+1}$.
\end{proof}

\subsection{Toeplitz subshifts and projections onto odometers}

Let $\Sigma$ be a compact Hausdorff space, $X = \Sigma^{\zz}$,
$S : X \rightarrow X$ is the left shift on $X$.

Assume $x = (x(n)) \in X$ is non--periodic recurrent point, $\{p_{i} \in
\EssP{x}\}$ is a sequence which complies with all conditions of
definition~\ref{per_struct_2}.

Let us consider a family of sets
\begin{equation}
\begin{split}
A^{i}_{j} & = \{ y(n) \in \Cl{\Orb(y)} \;|\;
    y(k+j) = x(k) \; \forall k \in \Per{p_{i}}(x) \} = \\
& = \{ y(n) \in \Cl{\Orb(y)} \;|\;
    y(n) = x(k) \; \forall n \equiv k+j \pmod {p_{i}} \,,\;
    k \in \Per{p_{i}}(x) \} \,, \\
& \quad j \in \{0, 1, \ldots, p_{i}-1\} \,,
    \quad i \in \nn \,.
\end{split}
\end{equation}
We can describe $A^{i}_{j}$ as the set of all points from $\Cl{\Orb (x)}$
which have the same $p_{i}$--skeleton with $S^{j}(x)$.

\begin{lemma}[compare with lemma 2.3 from~\cite{Williams}] \label{lemma_5}
The family of sets $\{A^{i}_{j}\}$ complies with the following properties
\begin{itemize}
    \item[(i)] For every $i \in \nn$ the family
        $\{ A^{i}_{j} \}_{j=0}^{p_{i}-1}$ is the periodic partition
        of the dynamic system $(\Cl{\Orb (x)}, S)$ of length $p_{i}$.
    \item[(ii)] $A^{i}_{n} \supset A^{j}_{m}$ when $i < j$ and
        $m \equiv n \pmod p_{i}$.
\end{itemize}
\end{lemma}

\begin{proof}
We mark first that for every $y \in X$ and for all $q_{1}$, $q_{2}$,
such that $\Per{q_{1}} (y)$, $\Per{q_{2}} (y) \neq \emptyset$,
the following implication is valid
\begin{equation} \label{eq_diez}
\left(
    q_{1} \text{ divides } q_{2}
\right)
\Rightarrow
\left(
    \Per{q_{1}} (y, \sigma) \subseteq \Per{q_{2}} (y, \sigma) \,, \;
    \sigma \in \Sigma
\right) \,.
\end{equation}

Let $x \in X$ and $\{ p_{i} \in \EssP{x} \}_{i \in \nn}$ satisfy
to requirements of lemma. Then from theorem~\ref{th_1} it follows
that the dynamic system $(\Cl{\Orb (x)}, S)$ is Toeplitz
and specially it is minimal.

We fix $i \in \nn$.

From lemma~\ref{lemma_3} it immediately follows that
$$
\Cl{\Orb (x)} = \bigcup_{j=0}^{p_{i}-1} A^{i}_{j}
$$
and the family of sets $\{ A^{i}_{j} \}_{j=0}^{p_{i}-1}$ satisfies
to the requirement~(iv) of definition~\ref{defn_9}.

Verify now the validity of requirement~(iii) of this definition. Assume that
$A^{i}_{j} \cap A^{i}_{k} \neq \emptyset$ for some $j \neq k$. Then
from lemma~\ref{lemma_3} and definition of the set $\Per{p_{i}} (x, \sigma)$
we get $x(n) = y(n+j) = y(n+k)$ for all $n \in \Per{p_{i}} (x)$ and
$$
\Per{p_{i}} (x, \sigma) = \Per{p_{i}} (y, \sigma) - j =
\Per{p_{i}} (y, \sigma) - k \quad \forall \sigma \in \Sigma \,.
$$
From corollary~\ref{corr_2} we have $p_{i} \in \EssP{y}$. Hence,
$p_{i}$ divides $|j-k|$ by definition of essential period. And it
contradicts to the inequality $0 < |j-k | < p_{i}$.

Let us prove now property~(ii) of definition~\ref{defn_9}.

From definition of sets $A^{i}_{j}$ the relations follow
\begin{equation} \label{eq_star_2}
\begin{split}
S(A^{i}_{j-1}) & \subseteq A^{i}_{j} \,,
    \quad j \in \{ 1, \ldots, p_{i}-1 \} \,; \\
S(A^{i}_{p_{i}-1}) & \subseteq A^{i}_{0} \,.
\end{split}
\end{equation}
With the help of these relations we immediately conclude that
\begin{equation} \label{eq_star_3}
S^{p_{i}} (A^{i}_{j}) \subseteq A^{i}_{j} \,, \quad
    j \in \{ 0, 1, \ldots, p_{i}-1 \} \,.
\end{equation}
The map $S$ is a homeomorphism. Hence, if even at least one
of the inclusions~(\ref{eq_star_2}) is strict, then
$$
S^{p_{i}} (A^{i}_{j}) \subsetneqq A^{i}_{j} \,, \quad
    j \in \{ 0, 1, \ldots, p_{i}-1 \} \,.
$$

From this remark and property~(iii) of definition~\ref{defn_9}, which we
have already verified, we conclude that in this case
$$
S^{p_{i}} \left( \Cl{\Orb(x)} \right) =
S^{p_{i}} \left( \bigcup_{j=0}^{p_{i}-1} A^{i}_{j} \right) =
\bigcup_{j=0}^{p_{i}-1} S^{p_{i}} \left( A^{i}_{j} \right) \subsetneqq
\bigcup_{j=0}^{p_{i}-1} = \Cl{\Orb (x)} \,.
$$

Since the set $\Cl{\Orb (x)}$ is Hausdorff and compact and
$S^{p_{i}}$ is a homeomorphism, then
$$
K = \bigcap_{m \geq 0} S^{m p_{i}} \left( \Cl{\Orb (x)} \right) \neq
\emptyset
$$
is the proper closed invariant subset of the dynamic system
$(\Cl{\Orb (x)}, S)$ contrary to the minimality of it.

Consider now property~(i) of definition~\ref{defn_9}.

All sets $A^{i}_{j}$ are closed. Really, we fix
$j \in \{ 0, 1, \ldots, p_{i}-1 \}$ and a convergent sequence
$y_{k} = y_{k} (n) \in A^{i}_{j}$. Let $y \in \Cl{\Orb (x)}$ is a limit
of this sequence. Since we have $y_{k}(m) = x(m-j)$, $k \in \nn$
for all $m \in \Per{p_{i}} (x)+j$, then proposition~\ref{prop_1} guarantees
$$
y(m) = x(m-j) \quad \text{for } m \in \Per{p_{i}} (x) + j \,.
$$
Consequently, $y \in A^{i}_{j}$ and the sets $A^{i}_{j}$ are closed.
That is $\{ A^{i}_{j} \}_{j=0}^{p_{i}-1}$ is the closed finite partition
of the dynamic system $(\Cl{\Orb (x)}, S)$. Therefore, each set
$A^{i}_{j}$ is open--closed in $\Cl{\Orb (x)}$ in the induced
topology.

The property~(ii) of lemma immediately follows from definition of sets
$A^{i}_{j}$, relations~(\ref{eq_diez}) and~(\ref{eq_star_3}) and
from lemma~\ref{lemma_3}.
\end{proof}

Let an odometer $(A, g)$ is built with the help of the sequence
$\{ p_{i} \}$.

Assume
$$
\vec{a} = (n_{i}) \in \prod_{i \in \nn} \zz_{p_{i}} \,.
$$
We denote
$$
A_{\vec{a}} = \bigcap_{i \in \nn} A^{i}_{n_{i}} \,.
$$

From the condition~(ii) of lemma~\ref{lemma_5} it immediately follows
that
\begin{equation} \label{eq_15}
\left( A_{\vec{a}} \neq \emptyset \right)
\Leftrightarrow
\left( \vec{a} \in A \subset \prod_{i \in \nn} \zz_{p} \right) \,.
\end{equation}
The condition~(i) of lemma mentioned guarantees that the family of
sets $\{ A_{\vec{a}} \,,\; \vec{a} \in A \}$ is partition of the
space $\Cl{\Orb (x)}$ and
\begin{equation} \label{eq_16}
S(A_{\vec{a}}) = A_{\vec{a}+\vec{e}}
\end{equation}
for every $\vec{a} \in A$.

Consider the following correspondence
\begin{equation*}
\begin{split}
    \pi & : \Cl{\Orb (x)} \rightarrow A \,; \\
    \pi & : A_{\vec{a}} \mapsto \vec{a} \,, \quad \vec{a} \in A
        \,.
\end{split}
\end{equation*}

From correlation~(\ref{eq_15}) we consequence the correctness of
this definition and formula~(\ref{eq_16}) guarantees the equality
$\pi \circ S = g \circ \pi$.

Mark that the map $\pi$ is continuous since $\pi^{-1}(V_{x_{j}})
= A^{j}_{x_{j}}$ for all $n \in \nn$, $x_{j} \in \zz_{p_{j}}$. In
other words all sets from the family~(\ref{eq_3.3}), which as we
know forms the base of topology of the space $A$, have
open--closed preimages in $\Cl{\Orb (x)}$ according to
lemma~\ref{lemma_5}.

\begin{theorem}
Assume that a point $x \in X$ is recurrent and a sequence
$\{ p_{i} \}_{i \in \nn}$ is a periodic structure on $x$ in
sense of definition~\ref{per_struct_2}.

Then the odometer $(A, g)$ built with the help of the sequence
$\{ p_{i} \}$ is an almost one-to-one factor of the flow
$(\Cl{\Orb (x)}, S)$ under the mapping $\pi$.

Moreover, two following conditions are equivalent:
\begin{itemize}
    \item[1)] a sequence $y \in \Cl{\Orb (x)}$ is Toeplitz;
    \item[2)]  $\pi^{-1}(\pi)(y) = \{y\}$.
\end{itemize}
\end{theorem}

\begin{proof}
Theorem is proved similarly to theorem~2.2 from~\cite{Williams}
(the single change is that the above lemma~\ref{lemma_5} must be
referred to instead of lemma~2.3 from~\cite{Williams}).
\end{proof}

\end{document}